\documentclass[11pt,a4paper,final]{article}
\usepackage{amssymb}
\usepackage[T1]{fontenc}
\usepackage[utf8]{inputenc}
\usepackage[english]{babel}
\usepackage{amsfonts}
\usepackage{enumitem}
\usepackage{amsmath}
\usepackage{amsthm}
\usepackage{booktabs}
\usepackage[colorlinks=true,allcolors=blue]{hyperref}
\usepackage[normalem]{ulem}
\usepackage{rotating}
\usepackage{amssymb}
\usepackage{tabularx,ragged2e,booktabs,caption}
\usepackage{geometry}
\usepackage{multirow}
\usepackage{subcaption}
\usepackage{stmaryrd}
 \usepackage{graphicx}

\def\R{\mathbb{R}}

\def\E{\mathbb{E}}

\def\R{\mathbb{R}}

\def\Z{\mathbb{Z}}

\def\1{\mbox{I\hspace{-.6em}1}} 

\def\1{\mbox{\hspace{.2em}I\hspace{-.6em}1}} 

\def\limiteasn{\renewcommand{\arraystretch}{0.5}
\begin{array}[t]{c}\stackrel{a.s.}{\longrightarrow} \\
{\scriptstyle
n\rightarrow+\infty}\end{array}\renewcommand{\arraystretch}{1}}

\DeclareMathOperator{\argmin}{argmin}
\DeclareMathOperator{\argmax}{argmax}
\newtheorem{theo}{Theorem}[section]
\newtheorem{lemma}{Lemma}
\newtheorem{proposition}{Proposition}
\numberwithin{equation}{section}

\newtheorem{algo}{Algorithm}
\newtheorem{remark}{Remark}

\newcommand*\interior[1]{\overset{\mathsf{o}}{#1}}

\newcommand{\ii}{\textbf{\textnormal{I}}}
\newcommand{\oo}{\textbf{\textnormal{O}}}

\def\limiten{\renewcommand{\arraystretch}{0.5}
\begin{array}[t]{c}
\stackrel{}{\longrightarrow} \\
{\scriptstyle n \rightarrow\infty}
\end{array}\renewcommand{\arraystretch}{1}}

\def\limiteproban{\renewcommand{\arraystretch}{0.5}
\begin{array}[t]{c}
\stackrel{{\mathcal P}}{\longrightarrow} \\
{\scriptstyle n \rightarrow\infty}
\end{array}\renewcommand{\arraystretch}{1}}

\def\limitesur{\renewcommand{\arraystretch}{0.5}
\begin{array}[t]{c}
\stackrel{{a.s.}}{\longrightarrow} \\
{\scriptstyle n \rightarrow\infty}
\end{array}\renewcommand{\arraystretch}{1}}

\def\limitedis{\renewcommand{\arraystretch}{0.5}
\begin{array}[t]{c}
\stackrel{{\mathcal D}}{\longrightarrow} \\
{\scriptstyle n \rightarrow\infty}
\end{array}\renewcommand{\arraystretch}{1}}

\begin{document}
 \title{General Hannan and Quinn Criterion for Common Time Series}
 \author{ BY Kare KAMILA\footnote{ This author has received funding from the European Union's Horizon 2020 research
 		and innovation programme under the Marie Sklodowska-Curie grant agreement No 754362.}}
 \maketitle
\vspace{0.25cm }
\hspace{1.7cm }\textit{SAMM, Université Paris 1 Panthéon-Sorbonne, FRANCE}
\vspace{0.1cm}
%
%
\begin{abstract}
This paper aims to study  data driven model selection criteria for a large class of time series, which includes  ARMA or AR($\infty$) processes, as well as
GARCH or ARCH($\infty$), APARCH and many others processes.  We tackled the challenging issue of designing adaptive criteria which enjoys the strong consistency property. When the observations are generated from one of the aforementioned  models, the new criteria, select the true model almost surely asymptotically. The proposed  criteria are based on the minimization of a penalized contrast akin to the Hannan and Quinn's criterion and then involved a term which is known for most classical time series models and for more complex models, this term can be data driven calibrated. 
Monte-Carlo experiments and an illustrative example on the CAC 40 index are performed to highlight the obtained results.
 \end{abstract}
\hspace{0.92cm} \textbf{Key words}: Time series, Model selection, consistency, data driven, HQ criterion.
\section{\textsc{Introduction}}

A common solution in model selection is to choose the model, minimizing a penalized based criterion which is the sum of two terms: the first one is the empirical risk (least squares, likelihood) that measures the goodness of fit and the second one is an increasing function of the complexity which aims to penalize large models and control the bias.
~\\
Therefore a challenging task when designing a penalized criterion is the specification of the penalty term. Considering leading model selection criteria (BIC, AIC, Cp, HQ to name a few), one can see that the penalty term is a product of the model dimension with a sequence which is specific  to the criteria. Indeed, a criterion is designed according to the goal one would like to achieve. The classical properties for model selection criteria include \textit{consistency},  \textit{efficiency} (oracle inequality, asymptotic optimality),  \textit{adaptative in the minimax sense}. 
\\

In this paper, we focus on consistency property which aims at identifying the data generating process with high probability or almost surely.  Hence, it requires the assumption whereby there exists a true model in the set of competitive models and the goal is to select this with probability approaches one as the sample size tends to infinity. In \cite{bar2019}, they studied model selection criteria regarding consistency in a  large class of time series, which is the interest of this paper. The leading criterion obtained in this framework is the BIC; with a relatively heavy penalty, it ensures the selection of quite simple models. Moreover, several papers have established the consistency property in particular settings. For instance, \cite{hannan} shows that the Hannan and Quinn (HQ) penalty $c\, \log \log n$ with $c>2$ leads to a consistent choice of the true order in the framework of AR type models. One year later, \cite{Hannan1980} (or \cite{han}) extended this result for  ARMA models. 

Also, it has been proven in several contexts, that  the BIC criterion \cite{schwarz} enjoys the consistency property: \cite{shibConsistency} in the density estimation using hypothesis testing for autoregressive moving average models, \cite{lebarbier} in density estimation for independent observations, \cite{bar2019} for a general class of time series, to name a few.


Compare to HQ penalty, the BIC penalty does not have the slowest rate of increase and then it can very often choose very simple models possible wrongs for small samples \cite{hannan}.
Moreover, the HQ criterion has been derived for linear time series:   AR models  in \cite{hannan}, ARMA models in \cite{Hannan1980} and \cite{han}. Is the HQ penalty still strongly consistent for heteroscedastic nonlinear models such as GARCH, APARCH or ARMA-GARCH? And what about a general class including linear and non linear models as well?

That raises a challenging question of designing robust penalties  for most classical time series models enjoying the model selection consistency. This is the  issue we want to address in this paper for a general class of times series called affine causal and defined below.
\\
~\\
\textbf{Class} $\mathcal{AC}(M,f):$ A process $X=(X_t)_{t\in \Z}$ belongs to $\mathcal{AC}(M,f)$ if it satisfies:
\begin{equation}
X_t=M\big((X_{t-i})_{i\in \mathbb{N}^*}\big)\, \xi_t+f\big((X_{t-i})_{i\in \mathbb{N}^*}\big) \;\; \mbox{for any}~ t\in \Z.
\label{eq:serie2}
\end{equation}
where $(\xi_t)_{t\in T}$ is a sequence  of zero-mean independent identically distributed random vectors (i.i.d.r.v) satisfying $\mathbb{E}(|\xi_0|^r)<\infty $ with $r\geq 1$  and $M$, $f$ : $\R^\infty \to \R$ are two measurable functions. \\
~\\
\noindent For instance,
\begin{itemize}
	\item if $M\big((X_{t-i})_{i\in \mathbb{N}^*}\big)=\sigma$ and $f\big((X_{t-i})_{i\in \mathbb{N}^*}\big)=\sum_{i=1}^{\infty} \phi_i X_{t-i}$,
	then $(X_t)_{t \in \Z}$ is an AR$(\infty)$ process;
	\item if $M\big((X_{t-i})_{i\in \mathbb{N}^*}\big)=\sqrt{a_0+a_1X^2_{t-1}+\cdots+a_pX^2_{t-p}}$ and $f\big((X_{t-i})_{i\in \mathbb{N}^*}\big)=0$,
	then $(X_t)_{t \in \Z}$ is an ARCH$(p)$ process.
\end{itemize}

\noindent Note that, numerous classical  time series models such as  ARMA($p,q$), GARCH($p,q$), ARMA($p,q$)-GARCH($p,q$) (see \cite{ding} and \cite{ling})
or APARCH$(\delta,p, q)$ processes (see \cite{ding}) belongs to $\mathcal{AC}(M,f)$.

\noindent The study of this type of process more often requires the classical regularity conditions on the functions $M$ and $f$, which are not restrictive at all and remain valid in various time serie models. Let us recall these conditions for $\Psi_{\theta}=f_{\theta}$ or  $M_{\theta}$ and $\Theta $ a compact set. 
\medskip

\noindent \textbf{Hypothesis A}$(\Psi_{\theta},\Theta)$: {\it Assume that $\|\Psi_{\theta}(0)\|_\Theta < \infty$ and there exists a sequence of non-negative real numbers $\big(\alpha_k(\Psi_{\theta},\Theta)\big)_{k\ge 1}$ such that $\sum_{k=1}^{\infty}\alpha_k(\Psi_{\theta},\Theta)< \infty$ satisfying:
	\[
	\|\Psi_{\theta}(x)-\Psi_{\theta}(y)\|_\Theta \le \sum_{k=1}^{\infty}\alpha_k(\Psi_{\theta},\Theta)|x_k-y_k| \; for \; all\; x,y \in \mathbb{R}^\mathbb{\infty}.
	\]}

\medskip

\noindent In addition, if the noise $\xi_0$ admits $r$-order moments  (for $r\geq 1$), let us define:
\begin{multline}
\Theta(r)=\Big \{\theta  \in \R^d,~ A(f_{\theta},\{\theta\})\; \textnormal{and}\; A(M_{\theta},\{\theta\})\; \textnormal{hold with} \\
\sum_{k=1}^{\infty} \alpha_k(f_{\theta},\{\theta\}) +\|\xi_0\|_r \, \sum_{k=1}^{\infty} \alpha_k(M_{\theta},\{\theta\}) < 1 \Big \}.
\end{multline}
Under this assumption, \cite{dou} showed that  there exists a stationary and ergodic solution to (\ref{eq:serie2}) with $r$-order moment for any $\theta \in \Theta(r)$. Moreover, \cite{barW} studied the consistency and the asymptotic normality of the QMLE of $\theta^*$ for  $\mathcal{AC}(M_{\theta^*},f_{\theta^*})$ . \\
~\\
\noindent The main contribution of this paper is the generalization of the HQ criterion to affine causal class: we provide a  minimal multiplicative penalty term $c_{min}$ so that all penalties of the form $2\,c\, \log \log n\, D_m$ with  $c\ge  c_{min}$ ensure the strong consistency property for affine causal models under some mild conditions on the Lipschitz coefficients of functions $M_{\theta},f_{\theta}$ ($D_m$ denotes the size of the model $m$). Monte Carlo experiments have been conducted in order to attest the accuracy of our new criteria. 
  
 
\medskip
\noindent The paper is organized as follows.
The model selection consistency along with notations and assumptions are described in Section
\ref{sec:cons}. Numerical results are presented in Section \ref{sec:monte} and Section \ref{sec:proof} contains the proofs.

\section{\textsc{Model Selection Consistency}}\label{sec:cons}
\subsection{Model Selection Procedure}

Let assume  $(X_1,\ldots,X_n)$ be a trajectory of a stationary affine causal process $m^*:=\mathcal{AC}(M_{\theta^*},f_{\theta^*})$, where $\theta^*$ is unknown. The goal of the consistency property is to come up with this true model  given a set of candidate model $\cal M$ such that $m^* \in \cal M$. 
\\
~\\
A $D_m$-dimensionnal model $m \in \cal M$ can be viewed as a set of causal functions $(M_{\theta},f_{\theta})$ with $\theta \in \Theta(m) \subset \R^{D_m}$. $\Theta(m) $ is the parameter set of the model $m$.


\noindent The consistency property will be study using the MLE first. Extension to QMLE will be done afterwards.
~\\
The MLE is derived from the conditional (with respect to the filtration $\sigma \big \{(X_{t})_{t\leq 0} \big \}$) log-likelihood of $(X_1,\ldots,X_n)$ when $(\xi_t)$
is supposed to be a Gaussian standard white noise.
Due to the linearity of a causal affine process, we deduce that this conditional log-likelihood (up to an additional constant) $L_n$ is defined for all $\theta \in \Theta$ by:
\begin{equation}
L_n(\theta):=-\frac{1}{2}\sum_{t=1}^n q_t (\theta) ~ , ~ \textnormal{with} \; q_t (\theta):=\frac{(X_t - f_{\theta}^t)^2}{H_{\theta}^t} + \log(H_{\theta}^t)
\label{eq:eq1}
\end{equation}
where $f_{\theta}^t:=f_{\theta}(X_{t-1},X_{t-2},\cdots)$, $M_{\theta}^t:=M_{\theta}(X_{t-1},X_{t-2},\cdots)$ and $H_{\theta}^t=\big (M_{\theta}^t\big )^2$.

\medskip
From now on, we drop the Gaussian assumption of the noise. 
Let $\cal M$ a finite family of candidate models containing the true one $m^*$. According to Proposition 1 in \cite{bar2019}, all these models can be included into a big one with  parameter space $\Theta$.
\noindent For each specific model $m\in {\cal M}$, we define the Gaussian MLE   $\widehat{\theta}(m)$ as
\begin{equation}
\widehat{\theta}(m)= \underset{\theta \in \Theta(m)}{\argmax} \;L_n(\theta).
\label{eq:qmle}
\end{equation}
To select the true model $m\in {\cal M}$, we consider a penalized contrast $C(m)$ ensuring a trade-off between $-2$ times the maximized  log-likelihood, which decreases with the size of the model, and a penalty increasing with the size of the model. Therefore, the choice of the "best" model $\widehat{m}$ among the estimated can be performed by minimizing the following criteria
\begin{equation}
\widehat{m}=  \underset{m \in \mathcal{M}}{\argmin} \;C(m)\quad\textnormal{with}\quad
C(m)=-2\,L_n\big(\widehat{\theta}(m)\big)+ \kappa_n(m)
\label{eq:cri}
\end{equation}
where $(\kappa_n)_n$ an increasing sequence depending on the number of observations $n$ and the dimension $D_m$. There exist several possible choices of $\kappa_n(m)$ including
\begin{itemize}
	\item  $\kappa_n(m)= 2 c\,D_m\, \log \log n$ with $c>1$, we retrieve the HQ criterion \cite{hannan};
	\item  $\kappa_n=D_m\,\log n$, $C$ yields to BIC criterion \cite{schwarz};
	\item  $\kappa_n=2\,D_m$, $C$ is the AIC criterion \cite{akaike}.
\end{itemize}

Basically the principle is that by increasing the size, the likelihood increases also. The question is whether this increase in complexity is offset by a sufficient increase in likelihood. If the answer is no, then the least complex model is used, even if it is less likely. If the answer is affirmative, then we accept to work with a more complex model. Of course, all the difficulty lies in the choice of weights between likelihood and complexity, and thus ultimately in the specification of the penalty multiplicative term $\kappa_n$.

What is the better weighting term of the model complexity? The aim here is by leveraging the  increasing rate of the likelihood, to propose a data driven  $\kappa_n$ in order to guarantee the strong consistency property to our model selection procedure i.e.
\begin{equation}
\widehat{m} \limitesur m^*.
\label{eq:cons}
\end{equation}
%
\subsection{Assumptions}
Some mild conditions will be required to prove the consistency of the considered model selection criteria.\\
~\\
The following  assumption is well-known as the identifiability one and is always required in order to guarantee the unicity of the global maximum of the MLE at the true parameter $\theta^*$. That is:\\
~\\
\textbf{Assumption A1}: {\it  For all $\theta, \, \theta' \in \Theta_m$, $
	~(f_{\theta}^0=f_{\theta'}^0) \; and\;  (M_{\theta}^0=M_{\theta'}^0) \implies \theta=\theta'. $ }\\
~\\
Another required assumption concerns the differentiability of $\Psi_{\theta}=f_{\theta}$ or  $M_{\theta}$ on $\Theta$. This type of assumption has already been considered in order to apply the QMLE procedure (see \cite{barW}, \cite{strau}, \cite{white}). \\
~\\
The following condition provides the invertibility of the Fisher's information matrix of $(X_1,\ldots,X_n)$ and was used to prove the asymptotic normality of the QMLE (see \cite{barW}).\\
~\\
\textbf{Assumption A2}: {\it One of the families $(\partial f_{\theta}^t/\partial \theta^{(i)})_{1\le i\le D_{m^*}}$ or $(\partial H_{\theta}^t/\partial \theta^{(i)})_{1\le i\le D_{m^*}}$ is a.e. linearly independent.} \\
~\\
\noindent Note that the definitions of the conditional log-likelihood requires that their denominators do not vanish.
Hence, we will suppose in the sequel that the lower bound of $H_{\theta}(\cdot) = \big(M_{\theta}(\cdot) \big)^2$ (which is reached since $\Theta$ is compact) is strictly positive:
\medskip

\noindent \textbf{Assumption A3}: {\it $\exists \underline{h}>0$ such that $\underset{\theta \in \Theta}{\inf}(H_{\theta}(x)) \ge \underline{h}$ for all $x\in \mathbb{R}^{\infty}$.}

\medskip
\noindent Next we assume the existence of the eighth order moment of the noise.
\medskip

\noindent \textbf{Assumption A4}: $\E[\xi_0^8] < \infty$.

\medskip 

\noindent We end the list of assumptions by assuming a suitable relation between the Fisher Information matrix $G(\theta^*_m)$ and the limiting Hessian matrix of the log-likelihood  $F(\theta^*_m)$ defined as follows
$$\big(F(\theta^*_m)\big )_{i,j}=\mathbb{E}\Big[\frac{\partial^2 q_0(\theta^*_m)}{\partial \theta_i \partial \theta_j}\Big]  \quad \mbox{and}\quad  (G(\theta^*_m))_{i,j}=\mathbb{E}\Big[\frac{\partial q_0(\theta^*_m)}{\partial \theta_i} \frac{\partial q_0(\theta^*_m)}{\partial \theta_j} \Big],$$
with $\theta^*_m:=(\theta^*,0,\ldots,0)^\top \in \Theta(m)$.

\medskip
\noindent \textbf{Assumption A5}: There exist absolutes constants $\alpha_1$ and $\alpha_2$ such that for any $m \in \cal M$ verifying $m^* \subset m$, 
\begin{equation}
\textnormal{\textbf{1}}_m^\top \Sigma_{\theta^*_m} \textnormal{\textbf{1}}_m= \alpha_1\,D_m^1+\alpha_2\,D_m^2
\end{equation}
where $D_m^1$ and $D_m^2$ are two integers such that $D_m^1+D_m^2=D_m$, $\;\;$ $\textnormal{\textbf{1}}_m:=(1,1,\ldots,1)^\top \in \R^{D_m}$,
$\Sigma_{\theta^*_m}:=G(\theta^*_m)^{1/2} F(\theta^*_m)^{-1} G(\theta^*_m)^{1/2}$.
~\\

For most classical affine causal models, \textbf{A5} is verified (see Proposition \ref{dodo2}). However, for more complex models such as ARMA-GARCH  with $\mu_4 \ne 3$, $\Sigma_{\theta^*_m}$ is hard to handle.
\subsection{Consistency Result}
Before stating the main result of this section, we give  important intermediate results. All proof of the results stated in this subsection can be found in Section \ref{sec:proof}.

The following Proposition suggests the existence of a term that will be the keystone of this work. 

\begin{proposition}\label{dodo}
	Let $m^*$ any affine causal model. For any model $m$ with $\theta^*_m \in \Theta(m)$, and under \textbf{A1-A5}, there exist $\alpha_1$, $\alpha_2$, $D_m^1$, $D_m^2$ such that  
	\begin{equation}\label{eq:}
	\limsup\limits_{n \rightarrow \infty} \, \frac{L_n\big(\widehat{\theta}(m)\big)-L_n(\theta^*_m)}{2\,\log \log n} = \frac{1}{4}\big(\alpha_1\,D_m^1+\alpha_2\,D_m^2\big) \quad a.s.
	\end{equation}
\end{proposition}
~\\
For every $m \in \cal M$, let us denote by $c_{min}(m)$ the following term that will be used several times
\begin{equation}\label{def}
c_{min}(m):= \frac{1}{4}\big(\alpha_1\,D_m^1+\alpha_2\,D_m^2\big)
\end{equation}
Now we state a result which provides the values of both $\alpha_1$ and $\alpha_2$  for most classical affine causal models.

\begin{proposition}\label{dodo2}
	Under the assumptions and notation of Proposition \ref{dodo}, we have
	\begin{itemize} 
		\item If $\mu_4=\E[\xi_0^4]=3$ (for instance for Gaussian  noise), then $\alpha_1=2$,  $\alpha_2=2$ and $c_{min}(m)=\frac{1}{2}\,D_m$;
		\item If the parameter $\theta$ identifying an affine causal model $ X_t=M_{\theta}^t\, \xi_t+f_{\theta}^t$ can be decomposed as $\theta=(\theta_1, \theta_2)'$ with  $f_{\theta}^t=\widetilde{f}^t_{{\theta}_1}$ and $M_{\theta}^t=\widetilde{M}^t_{{\theta}_2}$, then $\alpha_1=2$,  $\alpha_2=\mu_4-1$ and
		$$c_{min}(m)=\frac{ 1}{2}\,D_m^1+\frac{\mu_4-1}{4}\,D_m^2$$
	\end{itemize}
\end{proposition}
\noindent The second configuration in Proposition \ref{dodo2} includes classical time series
\begin{itemize}
	\item GARCH($p,q$), APARCH$(\delta,p,q)$ type models and related ones, $c_{min}(m)=\frac{\mu_4-1}{4} D_m$;
	\item ARMA($p,q$) models, $c_{min}(m)=\frac{D_m}{2}$ if the variance of the noise is known and  $c_{min}(m)=\frac{D_m-1}{2}+ \frac{\mu_4-1}{4}  $ otherwise. 
\end{itemize}

~\\
We can now state the first main result of this paper. 

\begin{theo}\label{theo:cons}
	Let $(X_1,\ldots, X_n)$ be an observed trajectory of an affine causal process $X$ belonging to $\mathcal{AC}(M_{\theta^*},f_{\theta^*})$  where $\theta^*$ is an unknown vector belonging to $\Theta(r) \subset \R^{D_{m^*}}$. Let also $\cal M$ be a finite family of candidate models such that  $m^* \in \cal M$. If assumptions \textbf{A1-A5} hold, there exist  $\alpha_1$, $\alpha_2$, and a minimal constant $c_{min}:=\max\big(\frac{\alpha_1}{4},\frac{\alpha_2}{4}\big)$ such that
	\\
	~\\ 
	\noindent  for any $ \kappa_n(m)=2\,c\, D_m\, \log \log n$ with
	\begin{equation}\label{eq:pena}
	c\ge c_{min}
	\end{equation}
	it holds for the selected model $\widehat{m}$ according to \eqref{eq:cri}
	\begin{equation}
	\widehat{m}\limitesur m^*,
	\label{eq:F}
	\end{equation}
\end{theo}

~\\
\begin{remark}
	\begin{enumerate}
		\item For classical configurations as seen in Proposition \ref{dodo2}, this result gives a generalization of  Hannan and Quinn criterion.
		
		\item For more complex models, the values of $\alpha_1$ and $\alpha_2$ are unknowns (at least until a better relationship between matrix $F(\theta^*_m)$ and $G(\theta^*_m)$ is found) and  so $c_{min}$ is also unknown. In these cases, we propose to use adaptive methods such as slope heuristic algorithm or  dimension jump \cite{arlot} to calibrate $c_{min}$.
	\end{enumerate}
\end{remark}







\vspace{0.5cm}
\noindent Let us mention that our result  generalizes the strong consistency obtained by \cite{hannan} for AR models as the affine causal class also contains GARCH models. It furthermore generalizes the result \cite{han} for ARMA models.

~\\
Theorem \ref{theo:cons} gives a theoretical guarantee on the consistency of the model selection procedure. However, it does not say anything about the convergence (and its rate) of the parameter estimate resulting from the model selection $\widehat{m}$. The following results shows that the final estimate $\widehat{\theta}_{\widehat{m}}$ is consistent and verifies a CLT.

\begin{theo}\label{theo:2}
	Under the assumptions of Theorem \ref{theo:cons}, it holds
	\begin{equation}
	\widehat\theta(\widehat{m}) \limiteproban \theta^*.
	\label{eq:F3}
	\end{equation}

	Moreover
	\begin{eqnarray}
	\sqrt{n} \, \Big (\big (\widehat\theta(\widehat{m})\big )_i-(\theta^*)_i\Big )_{i \in m^*}\limitedis {\cal
		N}_{|m^*|}\big (0 \ , \Sigma_{\theta^*, m^*}\big )
	\label{eq:con2}
	\end{eqnarray}
	with $\Sigma_{\theta^*, m^*}:=\ F(\theta^*)^{-1} G(\theta^*)
	F(\theta^*)^{-1}$
\end{theo}
\vspace{1cm}
In this subsection, we have used the MLE contrast without any distribution assumption on the noise to derive a consistency property. However, the contrast $L_n$ as in \eqref{eq:eq1} depends on all the past values of the process $X$, which are unobserved. In the next subsection, we will propose an extension of Theorem \ref{theo:cons} based on QMLE which does not require knowledge of the initial values of the process.

\subsection{
	Extension of Theorem \ref{theo:cons} when using QMLE.}

The goal of this subsection is to sharpen the conditions on the sequence $\kappa_n$ found in \cite{bar2019}.
Before stating the result, let recall a little bit some definitions and notations used in \cite{bar2019} about QMLE. 
Following the definition of $L_n$, the quasi-likelihood $\widehat{L}_n$ is

\begin{equation}
\widehat{L}_n(\theta):=-\frac{1}{2}\sum_{t=1}^n \widehat{q}_t (\theta) ~ , ~ \textnormal{with} \; \widehat{q}_t (\theta):=\frac{(X_t - \widehat{f}_{\theta}^t)^2}{\widehat{H}_{\theta}^t} + \log(\widehat{H}_{\theta}^t)
\label{eq:eql}
\end{equation}
where $\widehat{f}_{\theta}^t:=\widehat{f}_{\theta}(X_{t-1},X_{t-2},\ldots)$, $\widehat{M}_{\theta}^t:=\widehat{M}_{\theta}(X_{t-1},X_{t-2},\ldots)$ and $\widehat{H}_{\theta}^t=\big (\widehat{M}_{\theta}^t\big )^2$.

\noindent For each specific model $m\in {\cal M}$, we define the Gaussian QMLE   $\widehat{\theta}(m)$ as
\begin{equation}
\widetilde{\theta}(m)= \underset{\theta \in \Theta(m)}{\argmax} \;\widehat{L}_n(\theta).
\label{eq:qmle2}
\end{equation}
The choice of the "best" model $\widehat{m}$ among the estimated can be performed by minimizing the penalized contrast $\widehat{C}(m)$
\begin{equation}
\widehat{m}=  \underset{m \in \mathcal{M}}{\argmin} \;\widehat{C}(m)\quad\textnormal{with}\quad
\widehat{C}(m)=-2\,\widehat{L}_n\big(\widetilde{\theta}(m)\big)+ \kappa_n(m).
\label{eq:cri2}
\end{equation}

In this framework, we do not consider long memory process  and then we define the class $\mathcal{H}(M_{\theta},f_{\theta})$ a subset of $\mathcal{AC}(M_{\theta},f_{\theta})$ in which  every process has Lispchitz coefficients satisfying the following conditions

$$\alpha_j (f_{\theta})+ \alpha_j (M_{\theta})+\alpha_j (\partial _\theta f_{\theta})+\alpha_j (\partial _\theta M_{\theta})=O(j^{-\gamma}) \quad \mbox{with} \quad \gamma> 2.$$

\medskip
It is then straightforward to see that every process in the class $\mathcal{H}(M_{\theta},f_{\theta})$ verifies the following condition:
\\
~\\
\noindent \textbf{Condition $\boldsymbol{K}(\Theta)$}:
\[
\sum_{k\ge e}\frac{1}{\log \log k}\,\sum_{j \ge k} \alpha_j (f_\theta,\Theta)+ \alpha_j (M_\theta,\Theta)+\alpha_j (\partial _\theta f_\theta,\Theta)+\alpha_j (\partial _\theta M_\theta,\Theta) < \infty.
\] 

This finding allows to propose a sharpen generalization of both Theorem 3.1 in our previous paper \cite{bar2019}  and a similar result in \cite{kengne}.
\\
~\\
Before stating the result using QMLE, it is important  as in \cite{bar2019} or \cite{barW}, to distinguish the special case of NLARCH$(\infty)$ processes which includes for instance GARCH$(p,q)$ or ARCH$(\infty)$ processes. In such case, let us define the class:\\
~\\
\textbf{Class} $\widetilde{\mathcal{AC}}(\widetilde H_\theta)$: A process $X=(X_t)_{t\in \Z}$ belongs to $\widetilde{\mathcal{AC}}(\widetilde H_\theta)$ if it satisfies:
\begin{equation}
X_t=\xi_t  \, \sqrt {\widetilde H_\theta \big((X^2_{t-i})_{i\in \mathbb{N}^*}\big)}\;\; \mbox{for any}~ t\in \Z.
\label{eq:serie3}
\end{equation}
Therefore, if $ M^2_\theta \big((X_{t-i})_{i\in \mathbb{N}^*}\big)=H_\theta \big((X_{t-i})_{i\in \mathbb{N}^*}\big)=\widetilde H_\theta \big((X^2_{t-i})_{i\in \mathbb{N}^*}\big)$
then, $\widetilde{\mathcal{AC}}(\widetilde H_\theta)={\mathcal{AC}}(M_\theta,0)$.
In case of the class $\widetilde{\mathcal{AC}}(\widetilde H_\theta)$, we will use the  assumption $A(\widetilde H_{\theta},\Theta)$.
The new set of stationary solutions is for $r\geq 2$:
\begin{equation}
\widetilde {\Theta}(r)=\Big \{\theta  \in \R^d,~ {A}(\widetilde H_{\theta},\{\theta\})~ \mbox{holds with}~
\big (\|\xi_0\|_r \big )^2 \, \sum_{k=1}^{\infty} \alpha_k(\widetilde H_{\theta},\{\theta\}) < 1 \Big \}.
\end{equation}
\\
~\\
Finally, we propose to restrict class $\widetilde{\mathcal{AC}}(\widetilde H_\theta)$ to  $\widetilde{\mathcal{H}}(\widetilde H_\theta)$ as done with $\mathcal{H}(M_{\theta},f_{\theta})$ by considering all the process checking the condition
$$\alpha_j (H_{\theta})+\alpha_j (\partial _\theta H_{\theta})=O(j^{-\gamma})  \quad \mbox{with} \quad \gamma>2 $$
so that 
\[
\sum_{k\ge e}\frac{1}{\log \log k}\,\sum_{j \ge k} \alpha_j (H_{\theta})+\alpha_j (\partial _\theta H_{\theta}) < \infty.
\] 
~\\
We can now state the second main result.
\begin{theo}\label{theo:cons2}
	Let $(X_1,\ldots, X_n)$ be an observed trajectory of an affine causal process $X$ belonging to $\mathcal{H}(M_{\theta^*},f_{\theta^*})$ (or $\widetilde {\mathcal{H}}(\widetilde H_{\theta^*})$) where $\theta^*$ is an unknown vector belonging to $\Theta(r) \subset \R^{D_{m^*}}$ (or $\widetilde  \Theta(r) \subset \R^{D_{m^*}}$). Let also $\cal M$ be a finite family of candidate models such that  $m^* \in \cal M$. If assumptions \textbf{A1-A5} hold, there exist  $\alpha_1$, $\alpha_2$, and a minimal constant $c_{min}:=\max\big(\frac{\alpha_1}{4},\frac{\alpha_2}{4}\big)$ such that
	\\
	~\\ 
	\noindent  for any $\kappa_n(m)=2\,c\, D_m\, \log \log n$ with
	\begin{equation}\label{eq:penal}
	c\ge c_{min}
	\end{equation}
	it holds for the selected model $\widehat{m}$ according to \eqref{eq:cri2}
	\begin{equation}
	\widehat{m}\limitesur m^*.
	\label{eq:F2}
	\end{equation}
\end{theo}

\vspace{0.7cm}
All the comments made about the Theorem \ref{theo:cons} remain valid here. Moreover, recently, \cite{kengne} requires heavy penalties to ensure the strong consistency for the process in the class $\mathcal{AC}(M_{\theta^*},f_{\theta^*})$. Indeed, according to \cite{kengne}, it is necessary that $\kappa_n$ verified $\kappa_n/\log \log n \limiten \infty$ to obtain \eqref{eq:F2} which is a stronger condition since the HQ  criterion does not fulfill this condition and it is well known that HQ criterion is strongly consistent (see for instance  \cite{hannan}). Moreover, the new penalties found in this paper does not satisfy this condition, yet we ensure strongly consistency.

Also, let us mention that for small samples, heavy penalties such as those in \cite{kengne} can very often choose very simple models possible wrongs  \cite{hannan}. 

\begin{remark}
	Our results show that, asymptotically heavy penalties such as BIC penalty will ensure the consistency property. However, in practise, these penalties are used for fixed $n$ and most for small samples and it is important to point out their drawbacks. Very often, in the family of competitive models $\cal M$, there are misspecified and underfitted models ($\mathcal{M}'$). Since the difference $-2 \,\widehat{L}_n(\widehat{\theta}_m)+2 \,\widehat{L}_n(\widehat{\theta}_{m^*})>0$ for every $m \in \mathcal{M}'$, making the penalty, heavy could offset this positivity and can lead to the selection of some underfitted models and then wrong models.  To be more convincing of that, see the simulations (DGP III) experiments in  Section \ref{sec:monte}.
\end{remark}

\subsection{Algorithm of Calibration of the minimal constant }\label{sous:1}
There  exist several ways to calibrate the minimal constant $c_{min}$ including the dimension jump (presented below) and the data-driven slope estimation. Indeed, once an estimation of $c_{min}$ is obtained, many studies advocates the choice of $2\, \widehat{c}_{min}$ which turn out to be optimal (\cite{massart2007}, \cite{arlot} among others). Now we present the dimension jump algorithm.

\begin{algo}\label{algo:1}
	\textbf{Dimension Jump \cite{arlot}}
	\begin{enumerate}
		\item Compute the selected model $\widehat m(c)$ as a function of $c > 0$
		$$\widehat m(c)=  \underset{m \in \mathcal{M}}{\argmin} \; \widehat{L}_n\big(\widehat{\theta}_m\big)+c\,\textnormal{pen}_{\mbox{shape}}(m); \quad \textnormal{pen}_{\mbox{shape}}(m)=D_m \,\log \log n$$
		\item Find $\widehat c_{min}$ such that $D_{\widehat m(c)}$ is "huge" for $c< \widehat{c}_{min}$ and "reasonnably small" for  $c\ge \widehat{c}_{min}$;
		\item Select the model $\widehat m:=\widehat m(2\, \widehat{c}_{min}).$
	\end{enumerate}
\end{algo}
This algorithm has been implemented in \cite{baudry} which gives several details including the grid for $c$ values.

Let us notice that, in view of obtaining penalties, there is no need to calibrate the $c_{min}$ constant for most classical time series models. However, since the fourth order moment of the noise is unknown, a consistent estimate of this term is required. To do that, we proceed as in the estimation of the variance of the noise as in the  Mallows Cp.

A consistent estimator $\widehat{\mu}_{\overline m,4}$ of  $\mu_4=\E[\xi_0^4]$ can be :
\[
\widehat{\mu}_{\overline m,4}:=\frac{1}{n} \sum_{t=1}^n \big(\widehat {\xi}_t(\overline m)\big)^4\quad \mbox{with}\quad \widehat {\xi}_t(\overline m):=\big (\widehat M^t_{\widehat \theta_{\overline m}}\big )^{-1}\big (X_t -\widehat f^t_{\widehat \theta_{\overline m}} \big )
\]
where we suppose that $\overline m$ is the "largest" model in the family ${\cal M}$, typically the largest order of a family of time series. As a result an estimator of the $c_{min}$ constant to consider in the penalty $\kappa_n$ is  
\begin{itemize}
	\item  $ \frac{\widehat{\mu}_{\overline m,4}-1}{4}$ for GARCH family and related ones;
	\item  $ \frac{1}{2}$ for ARMA family with known variance.
\end{itemize}

\section{\textsc{Numerical Experiments}}\label{sec:monte}
In this section, several numerical experiments are conducted to assess the consistency property (Section \ref{sec:cons}) of our new criteria.  
\subsection{Monte Carlo: Consistency}
This subsection studies the performance of the model selection criteria found in Section \ref{sec:cons}. We have considered three different Data Generating Process (DGP):
\begin{eqnarray*}
	\mbox{DGP I} &  & X_t=0.5\,X_{t-1}+0.2\,X_{t-2}+\xi_t,\\
	\mbox{DGP II} &  & X_t=\big(0.2+0.4\,X_{t-1}^2+0.2\,X_{t-2}^2\big)^{1/2} \, \xi_t,\\
	\mbox{DGP III} &  & X_t=0.1(\,X_{t-1}+X_{t-2}+\ldots+X_{t-6})+\xi_t,
\end{eqnarray*}
where $(\xi_t)$ will be a white Gaussian noise with variance one at first and a Student with 5 degrees of freedom on the other hand. 
\noindent For the first and the second DGP, we considered as competitive models all the models in the family ${\cal M}$ defined by:
$$
{\cal M}=\big \{\mbox{ARMA$(p,q)$ or GARCH$(p',q')$ processes with $0\le p,q,p' \le 5$, $1\le q'\le 5$} \big \}.
$$
Therefore, there are $66$ candidate models as in \cite{bar2019}. The goal is to compare the ability of selecting the true model for $\kappa_n=\log n\, D_m$, $\kappa_n=2\,\widehat{c}_{min}\, \log \log n\,D_m$ ( in accordance with the condition \eqref{eq:penal} and $\kappa_n=\textbf{2}\times 2\,\widehat{c}_{min}\, \log \log n\,D_m$.  
Moreover, from Theorem \ref{theo:2} $\widehat{c}_{min}$ does not need to be estimated and worth one half for Gaussian noise. But for Student noise, $\widehat{c}_{min}=\max\big(\frac{1}{2},\frac{\widehat{\mu}_{\overline m,4}-1}{4}\big)$ for the DGP I and   $\widehat{c}_{min}=\frac{\widehat{\mu}_{\overline m,4}-1}{4}$ for DGP II. 
The Table \ref{tab:1} presents the results of the selection procedure. As we can notice, the three penalties  have a good consistency property.  Moreover, for  $n$ relatively small, the penalty $\kappa_n=2\,\widehat{c}_{min}\, \log \log n\,D_m$ is better than both others. For larger $n$,  $\kappa_n=\widehat{c}_{min}\, \log \log n$ is the best the penalty to consider. 


\begin{center}
	\captionof{table}{ Percentage of selected order  based on 500 independent replications depending on sample's length for the penalty $\widehat{c}_{min}$, $2\, \widehat{c}_{min}$ and  $\log n$, where and W, T, O refers to wrong, true and overfitted selection.} \label{tab:1}
	\small{
		\begin{tabular}{|l|l|l|l|l|l|}
			\toprule
			\multicolumn{1}{|l|}{}& \multicolumn{1}{|l|}{$n$ }& \multicolumn{1}{c}{$100$ } & \multicolumn{1}{c}{$500$ }& \multicolumn{1}{c}{$1000$} & \multicolumn{1}{c|}{$2000$ }\\
			& &  $\widehat{c}_{min}$\, $2\,\widehat{c}_{min}$ \, $\log n$   &  $\widehat{c}_{min}$\, $2\,\widehat{c}_{min}$ \, $\log n$  & $\widehat{c}_{min}$\, $2\,\widehat{c}_{min}$ \, $\log n$  & $\widehat{c}_{min}$\, $2\,\widehat{c}_{min}$ \, $\log n$ \\
			\midrule
			&W &$58.2 $  \hspace{0.35cm}$80.8$ \hspace{0.35cm} $94.6 $ &$23.6$ \hspace{0.35cm}  $ 26$ \hspace{0.35cm} $36.2$ &$18.8$\hspace{0.35cm}  $17.6$\hspace{0.35cm}  $17.8$  &$4.2$\hspace{0.5cm}  $ 7.4$\hspace{0.8cm}$ 7.8$\\
			\small{DGP I }  &T &$32.6$  \hspace{0.4cm}$19$ \hspace{0.7cm}$5.4$  &$64.2$\hspace{0.35cm} $73.4$ \hspace{0.35cm}$63.8$ &$71.6$  \hspace{0.3cm}  $82.2$\hspace{0.35cm} $82.0 $   &$89.4$\hspace{0.4cm}   $92.2$\hspace{0.4cm} $92.2$\\
			\small{Gaussian }& O &\,$9.2$ \hspace{0.45cm}$0.2 $\hspace{0.8cm}$0   $ &$12.2$\hspace{0.4cm}  $0.6$\hspace{0.5cm}  $0 $ &$9.6$\hspace{0.55cm} $0.2$\hspace{0.7cm} $0.2$  &$6.4$\hspace{0.5cm}  $0.5$\hspace{0.75cm} $0$ \\
			\hline
			\midrule
			& W &$70.0$\hspace{0.4cm}$84.8$  \hspace{0.35cm}$94.0$ &$17.4$ \hspace{0.3cm} $21.6$  \hspace{0.3cm}$35.4$ &$16.4$\hspace{0.45cm}$ 15.8$\hspace{0.45cm}$ 15.8$  &$3.6$\hspace{0.5cm} $3.6$ \hspace{0.75cm}$5$    \\
			\small{DGP II }& T &$29.8$ \hspace{0.35cm}$15.2$  \hspace{0.35cm}$6.0$ &$81.0$ \hspace{0.3cm} $78.4$ \hspace{0.2cm} $64.6$ &$83.4$\hspace{0.35cm} $84.2$\hspace{0.35cm} $84.2$ &$96.4$\hspace{0.35cm}  $96.4$\hspace{0.45cm} $95$ \\
			\small{Gaussian}& O &\,$0.2 $ \hspace{0.55cm}$0$ \hspace{0.75cm}$0$ & $1.6$\hspace{0.65cm} $0$\hspace{0.75cm} $0.$ &$0.2$\hspace{0.7cm}$0$\hspace{0.75cm} $0.0$  &$0$\hspace{0.8cm}  $0$ \hspace{0.85cm}$ 0$\\
			\hline
			\midrule
			& W &$70.1$ \hspace{0.35cm}$86.8$ \hspace{0.35cm}$97.5$ &$38.3$\hspace{0.35cm} $53.1$\hspace{0.35cm} $35.0$ &$16.0$\hspace{0.35cm} $19.5$\hspace{0.35cm} $16.3$ &$18$\hspace{0.55cm}  $22.0$\hspace{0.75cm} $18$   \\
			\small{DGP I } &T &$19.8$ \hspace{0.35cm}$13.2$ \hspace{0.45cm}$2.5$ &$59.7$\hspace{0.35cm} $46.9$\hspace{0.35cm} $65.9$ &$83.9$\hspace{0.35cm} $80.5$\hspace{0.35cm} $82.7$ &$81.4$ \hspace{0.35cm}$78.0$\hspace{0.55cm} $82.0$\\
			\small{Student }& O &\,$10.1$\,\,\,\,\,\, \,$0$ \hspace{0.8cm}$0.$ & \,$2.0$\,\,\,\,\,\,\,\,\,\, $0$\hspace{0.7cm} $0.1$ &$0.1$\hspace{0.65cm} $0$\hspace{0.75cm} $0$  & $0.6$\,\,\,\,\;\;\;\;  $0$\hspace{1cm} $0$\\
			\hline
			\midrule
			& W &$75.0$ \hspace{0.35cm}$84.2$\hspace{0.35cm}$88.8$ &$45.6$\hspace{0.35cm} $57.2$\hspace{0.35cm} $55.4$ &$25.6$\hspace{0.35cm} $32.6$\hspace{0.35cm} $26.2$ &$14$\hspace{0.5cm}  $16.0$\hspace{0.75cm} $13.0$   \\
			\small{DGP II } &T &$22.4$ \hspace{0.35cm}$15.0$ \hspace{0.15cm} $11.2$ &$49.0$\hspace{0.35cm} $41.2$\hspace{0.35cm} $44.6$ &$71.0$\hspace{0.35cm} $66.8$\hspace{0.35cm} $73.8$ &$85.0$\hspace{0.3cm} $84.0$\hspace{0.65cm} $87.0$\\
			\small{Student }& O &\,$2.6$ \hspace{0.5cm}$0.8$ \hspace{0.5cm}$0.0$ &$4.4$\hspace{0.55cm} $1.6$\hspace{0.55cm} $1.0$& $3.4$\hspace{0.55cm} $0.6$\hspace{0.55cm} $0.0$  & $1.0$\hspace{0.65cm}  $0$\hspace{1.2cm}$ 0$\\
			\bottomrule
	\end{tabular}}
\end{center}
\vspace{0.5cm}

~\\
\noindent For DGP III, as we want to exhibit the possible "non consistency" of BIC for small samples, we have considered as the competitive set, the hierarchical family of AR models $$\mathcal{M}'=\big \{AR(1),\ldots, AR(15)\big\}$$
for $n=100,200,400,500,1000,2000$. In Table \ref{tab:2} below the percentage of selected order based on 1000 independent replications are presented for the above three penalties.

\begin{center}
	\captionof{table}{ Percentage of selected order  based on 1000 independent replications depending on sample's length using penalty terms $\kappa_n^1=\widehat{c}_{min}$, $\kappa_n^2=2\, \widehat{c}_{min}$ and  $\kappa_n^3=\log n$, for DGP III.} \label{tab:2}
	\small{
		\begin{tabular}{|l|l|l|l|l|l|l|}
			\toprule
			\multicolumn{1}{|l|}{$n$ }& \multicolumn{1}{c}{$100$ } & \multicolumn{1}{c}{$200$ }& \multicolumn{1}{c}{$400$} &
			\multicolumn{1}{c}{$500$ }& \multicolumn{1}{c}{$1000$} &
			\multicolumn{1}{c|}{$2000$ }\\
			&$\kappa_n^1$\hspace{0.35cm}$\kappa_n^2$\hspace{0.2cm} $\kappa_n^3$   &  $\kappa_n^1$\; $\kappa_n^2$ \; $\kappa_n^3$ & $\kappa_n^1$\;\;\; $\kappa_n^2$ \; $\kappa_n^3$  & $\kappa_n^1$\;\;\; \;$\kappa_n^2$ \;\; $\kappa_n^3$ & $\kappa_n^1$\;\;\; $\kappa_n^2$ \, $\kappa_n^3$ & $\kappa_n^1$\;\;\; $\kappa_n^2$ \;\; $\kappa_n^3$\\
			\midrule
			$p<6$ &61.5   96 \,\,99 &46\, 83.5\, 97\, &29.5\, 59\, 87.5 &19.5\, 45.5\, 76.5 &3.0  \, 9.5\, 38.5& \,0\;\;\;\, 0.5\;\;\;5.0 \\
			$p=6$  &15\;\;\; 3.0 1.0 &25 \,13\;\,\; 2.5 &44 \,\;37.5\,12.5  &50.5\;\;49.5\,\;22.5&68\;\;\;86.5\;\;61 &72.5\;96\;\;94.5 \\
			$p>6$  &23.5 1.0\;\;0 &29\;\,3.5\;\;\;0.5 &26.5\,3.5\;\;\;0 &\;30\;\;\;\,5.0\;\;\;\;1.0 &29\;\;\;\,4.0\;\;\;0.5 &27.5\,3.5\,\;\;0.5  \\
			\bottomrule
	\end{tabular}}
\end{center}
~\\
These results invite us to be cautious when using the BIC for small sample sizes, whereas the proposed adaptive penalty is more robust, as it at least allow us to recover an overfitted model that is less harmful than a wrong model most often chosen by the BIC. 

\subsection{Real Data Analysis: financial time series}
\begin{figure}
	\centering
	\includegraphics[width=10cm,height=7cm]{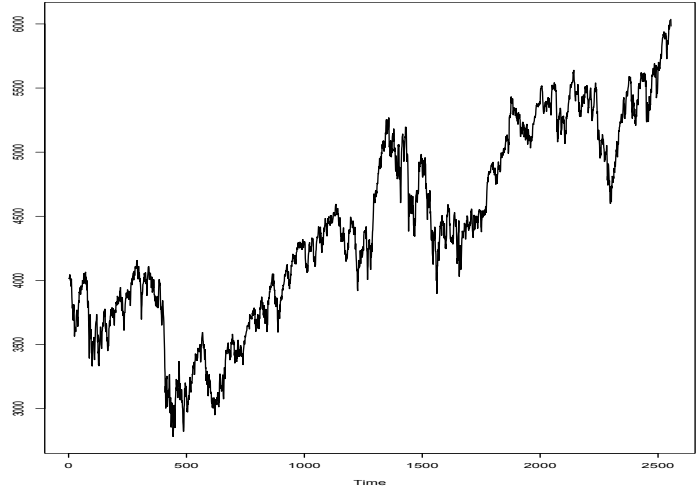}
	\caption{Daily closing CAC 40 index (January 1st, 2010 to December 31 st, 2019).}\label{fig:cac1}
\end{figure}

CAC 40 is a benchmark French stock market index and is highly regarded in many statistical studies . Let consider the  daily closing prices of the CAC 40 index from January 1st, 2010 to December 31st, 2019 plotted in Figure \ref{fig:cac1}. Over the period under review, the CAC 40 increased.

\begin{figure}[htb]
	\centering
	\begin{tabular}{cc}
		\subcaptionbox{Time plot of log returns.\label{fig:cac2}}[0.4\linewidth]{\includegraphics[width=7cm,height=7cm]{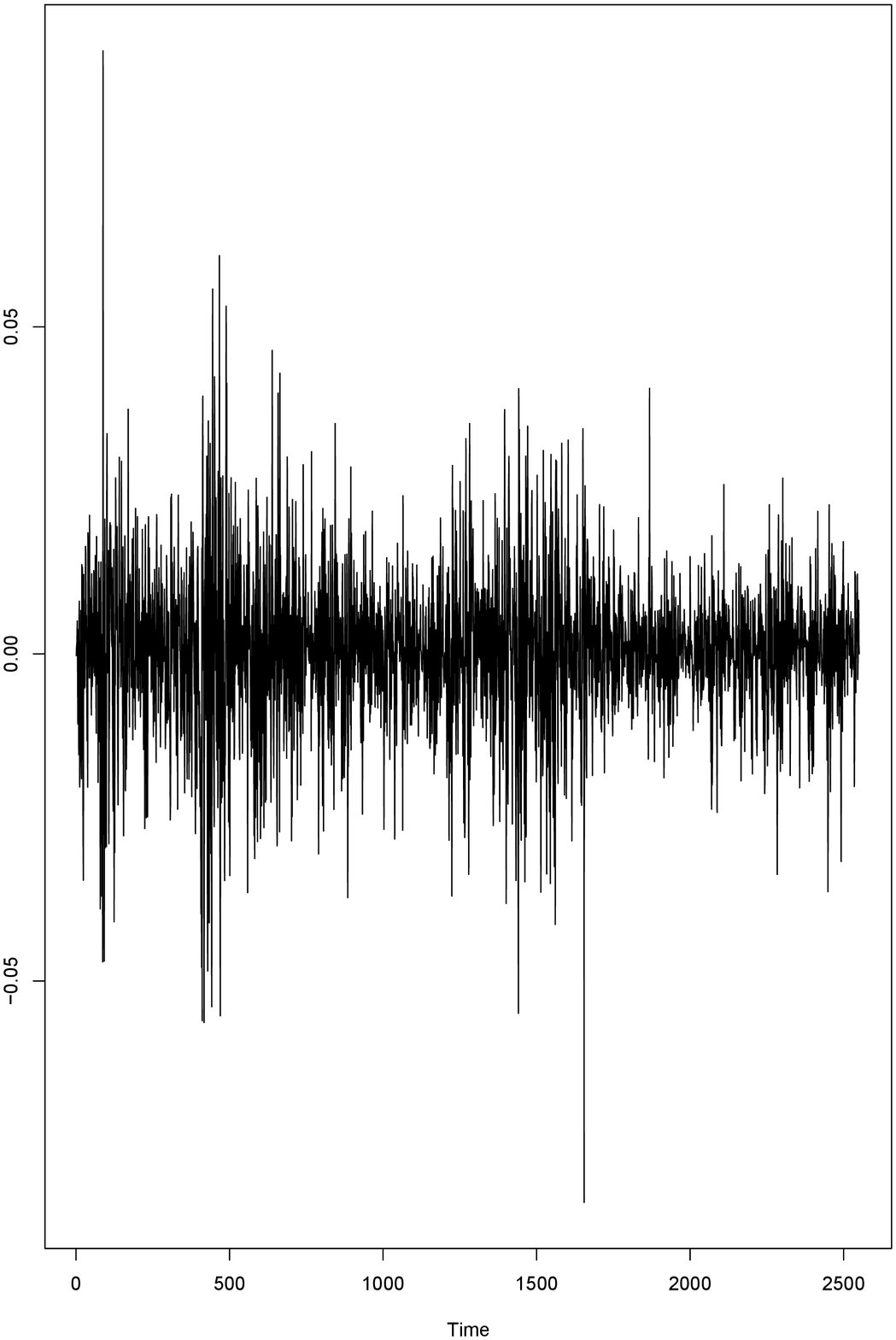}} &
		\subcaptionbox{Correlograms of log returns.\label{fig:cac3}}[0.4\linewidth]{\includegraphics[width=7cm,height=7cm]{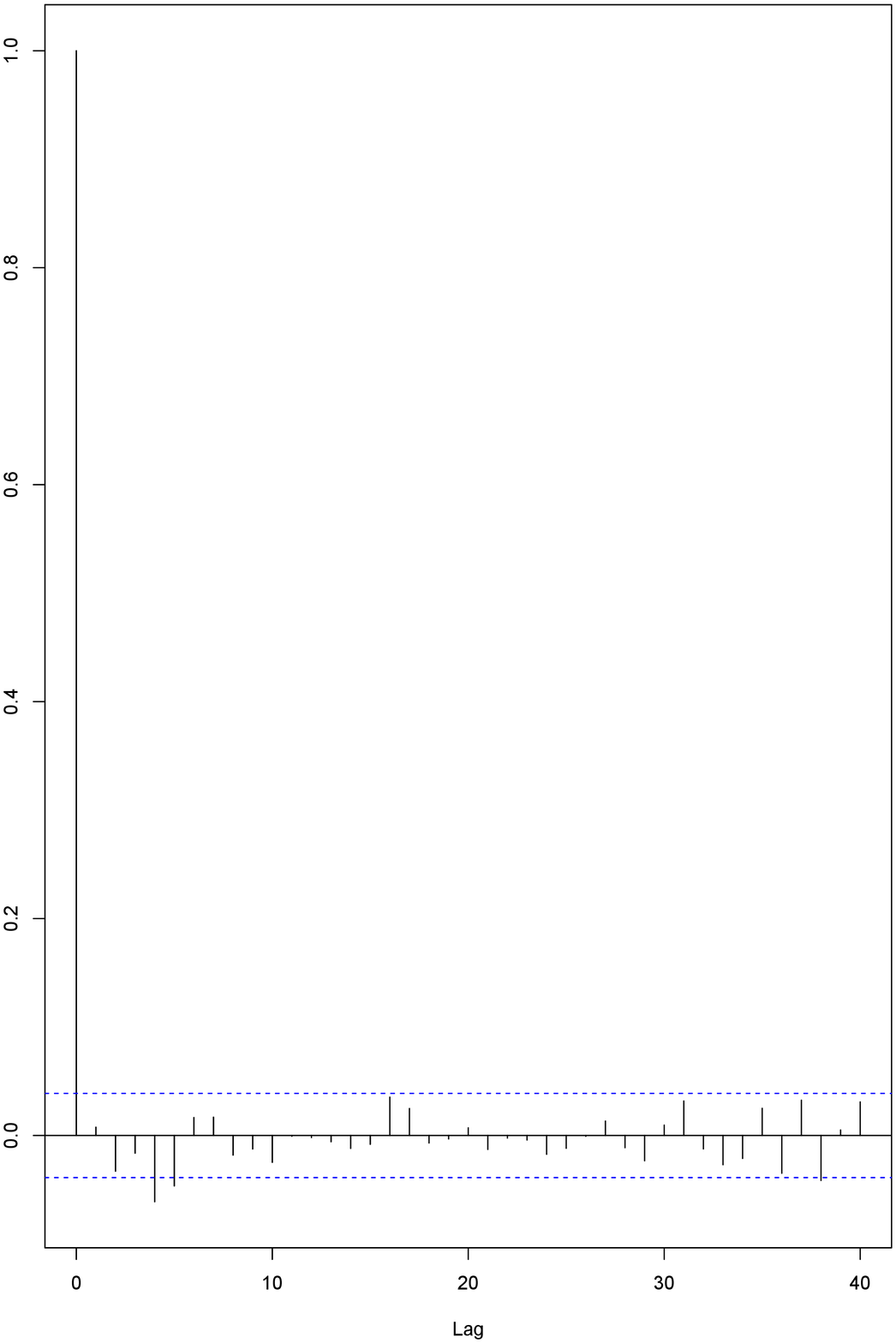}} 
	\end{tabular}
	\caption{Daily closing CAC 40 index (January 4th, 2010 to December 31 st, 2019).\label{arma}}
\end{figure}

\begin{figure}
	\centering
	\includegraphics[width=10cm,height=7cm]{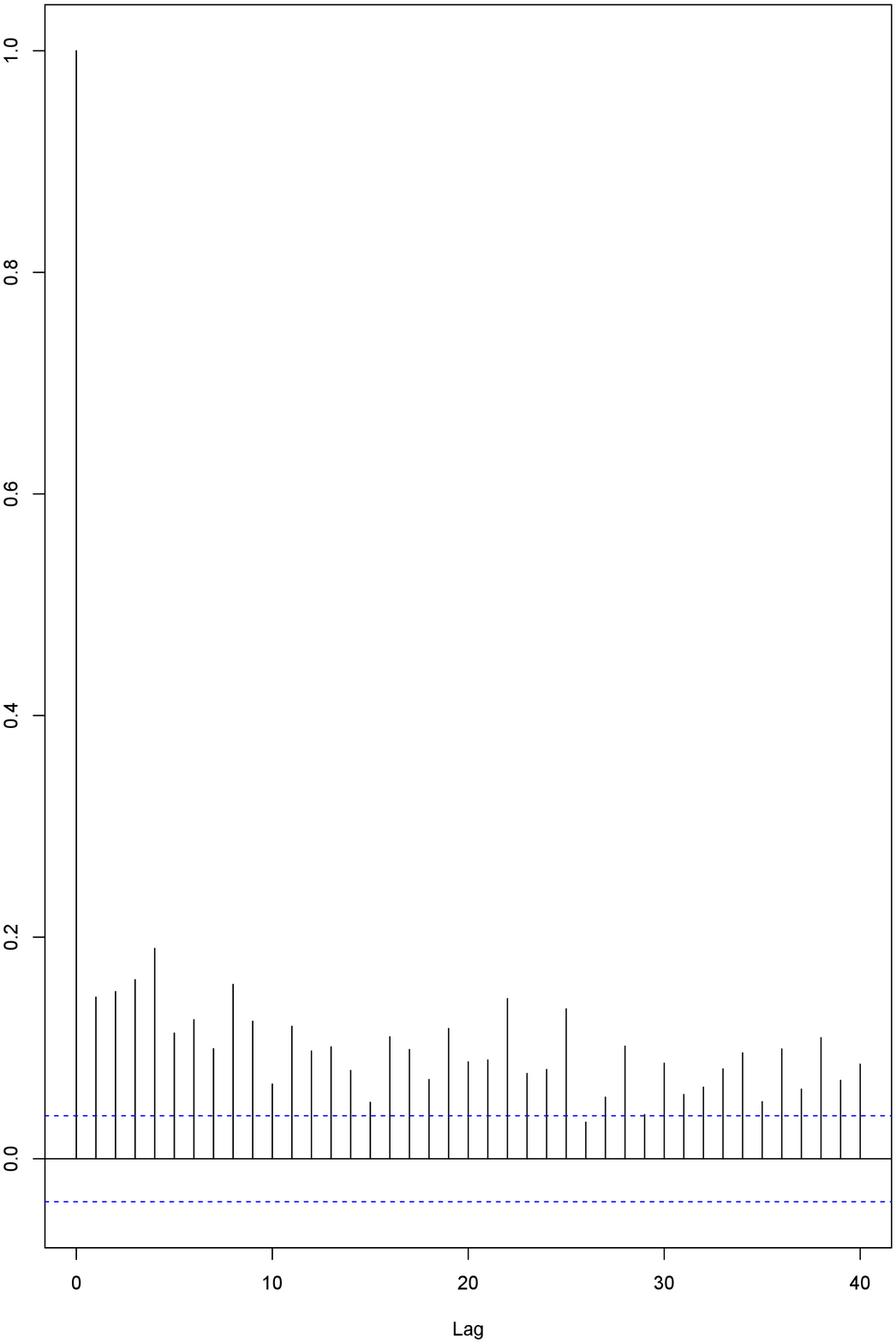} 
	\caption{Sample autocorrelations of squared returns of the CAC 40 index (January 1st, 2010 to December 31 st, 2019).}\label{fig:cac4}
\end{figure}

To analyze this type of data, it is common to consider the returns (see Figure \ref{arma}). We can see that the return values display some small auto-correlations. Also, from Figure \ref{fig:cac4}, the squared returns of CAC 40 are strongly auto-correlated. These facts suggest that the strong white noise assumption cannot
be sustained for this log-returns sequence of the CAC 40 index.

Hence, let consider the competitive set of models $\cal M$ used in the previous subsection in order to propose the best suitable model for these data:
$$
{\cal M}=\big \{\mbox{ARMA$(p,q)$ or GARCH$(p',q')$ processes with $0\le p,q,p' \le 5$, $1\le q'\le 5$} \big \}.
$$
Using the adaptive penalty and the BIC criterion, we find out that the GARCH($1,1$) is the best model over $\cal M$ with respect to both criteria. This fact is in accordance with \cite{francq2019} which found the same result using the returns of the CAC 40 index from March 2, 1990 to December 29, 2006.

\section{Proofs} \label{sec:proof}

\subsection{Proof of Theorem \ref{theo:cons}}
Usually, this proof is divided into two parts: one has to show that as $n$ tends to infinity the probability of overfitting goes to zero and so is the probability of misspecification.
\subsubsection{Overfitting Case}


Let  $m \in \cal M$  such as $m^* \subset m $. We want to show that $C(m^*) \le C(m)$ a.s. asymptotically. 
We have 
\begin{eqnarray}
\nonumber	C(m^*) \le C(m)  &\Longleftrightarrow& -2\,L_n\big(\widehat{\theta}(m^*)\big)+2c\,\log \log n \, D_{m^*} \le  -2\,L_n\big(\widehat{\theta}(m)\big)+2c\,\log \log n \, D_{m} \\
\label{rec}	&\Longleftrightarrow& \frac{L_n\big(\widehat{\theta}(m)\big)-L_n(\theta^*)}{2\,\log \log n} \le \frac{L_n\big(\widehat{\theta}(m^*)\big)-L_n(\theta^*)}{2\,\log \log n} + \;\;c\,(D_{m}-D_{m^*})
\end{eqnarray}
therefore, a necessary and sufficient  condition to avoid overfitting can be stated by taking $\limsup\limits_{n \rightarrow \infty}$ on both sides of the inequality \eqref{rec}; that is 
by virtue of definition \eqref{def}
\begin{equation}\label{eq:con}
c_{min}(m)-c_{min}(m^*) \le c\,(D_{m}-D_{m^*}) \quad \mbox{for} \quad m^*\subset m,
\end{equation}
i.e., $$  \frac{\alpha_1}{4}\,(D_m^1-D_{m^*}^1 )+\frac{\alpha_2}{4}\,(D_m^2-D_{m^*}^2)  \le c\,(D_{m}-D_{m^*})$$
which is fulfilled for any constant $c$ such as in \eqref{eq:pena}. Indeed, $c_{min}=\max\big(\frac{\alpha_1}{4},\frac{\alpha_2}{4}\big)$ and 
\begin{eqnarray*}
	\frac{\alpha_1}{4}\,(D_m^1-D_{m^*}^1 )+\frac{\alpha_2}{4}\,(D_m^2-D_{m^*}^2) & \le & c_{min}\Big (D_m^1-D_{m^*}^1 +D_m^2-D_{m^*}^2 \Big)\\
	& =& c_{min} (D_m-D_{m^*})
\end{eqnarray*}
where the inequality holds since $m^*\subset m$ that implies $D_m^1 \ge D_{m^*}^1 $ and $D_m^2 \ge D_{m^*}^2 $.
Hence the associated criterion  $\kappa_n$ will avoid overfitting.


\subsubsection{Misspecification/Underfitting Case} 
All misspecified/underfitted models are contained in the set  $$\mathcal{M}'=\Big\{m \in \mathcal{M} : (m^* \not \subseteq  m) \cup (m \subset m^*)  \Big\}.$$ The proof is exactly as the one done  in \cite{bar2019}. But for the sake of completeness, we give here some important steps of the proof.  

\noindent The goal is to show that for every $m \in \mathcal{M}'$, it holds 
\begin{equation}\label{eq:ous}
C(m^*) \quad < \quad C(m) \quad a.s.
\end{equation}
Let $m \in \mathcal{M}'$. From Proposition 2 in \cite{bar2019} and using continuous mapping Theorem 
\begin{equation}
\frac 1 n \Big[L_n\big(\widehat{\theta}(m^*)\big)-L_n\big(\widehat{\theta}(m)\big)\Big]= A(m_0)+ o_{a.s}(1)
\end{equation}
where 
$$A(m):=L(\theta^*)-L(\theta^*(m))  \qquad \mbox{with} \quad L(\theta)=-\frac 1 2 \,\E[q_0(\theta)].$$
Let us denote by $\mathcal{F}_t:=\sigma \big (X_{t-1},X_{t-2},\cdots \big )$. Using conditional expectation, we obtain
\begin{equation}\label{cond}
L(\theta^*)-L(\theta)=\frac{1}{2} \, \E \Big[ \E \big [q_0(\theta)-q_0(\theta^*)~|~\mathcal{F}_0 \big ]\Big].
\end{equation}
But,
\begin{eqnarray*}
	\E \big [q_0(\theta)-q_0(\theta^*)~|~\mathcal{F}_0 \big ]& =&
	\E\Big [\frac{(X_0 - f_{\theta}^0)^2}{H_{\theta}^0} + \log(H_{\theta}^0) -\frac{(X_0 - f_{\theta^*}^0)^2}{H_{\theta^*}^0} - \log(H_{\theta^*}^0) ~ |~\mathcal{F}_0 \Big ]\\
	&=& \log \Big (\frac{H_{\theta}^0}{H_{\theta^*}^0}\Big )+\frac{\E\big [(X_0 - f_{\theta}^0)^2~|~\mathcal{F}_0\big ]}{H_{\theta}^0} -\frac{\E\big [(X_0 - f_{\theta^*}^0)^2~|~\mathcal{F}_0\big ]}{H_{\theta^*}^0}\\
	&=&\log\Big (\frac{H_{\theta}^0}{H_{\theta^*}^0}\Big )-1+ \frac{\E\big [(X_0 -f_{\theta^*}^0+f_{\theta^*}^0- f_{\theta}^0)^2~ |~\mathcal{F}_0\big ]}{H_{\theta}^0} \\
	&=&\frac{H_{\theta^*}^0}{H_{\theta}^0}- \log\Big (\frac{H_{\theta^*}^0}{H_{\theta}^0}\Big )-1 + \frac{ (f_{\theta^*}^0-f_{\theta}^0)^2}{H_{\theta}^0}
\end{eqnarray*}
Thus from \eqref{cond},
\begin{eqnarray*}
	2\, A(m)& =& \E \Big [ \frac{H_{\theta^*}^0}{H_{\theta^*(m)}^0}- \log\Big (\frac{H_{\theta^*}^0}{H_{\theta^*(m)}^0}\Big )-1 + \frac{ (f_{\theta^*}^0-f_{\theta^*(m)}^0)^2}{H_{\theta^*(m)}^0} \Big ]\\
	& \ge & \E \Big [ \frac{H_{\theta^*}^0}{H_{\theta^*(m)}^0}\Big ]- \log\Big (\E \Big [ \frac{H_{\theta^*}^0}{H_{\theta^*(m)}^0}\Big ]\Big )-1 +\E \Big [  \frac{ (f_{\theta^*}^0-f_{\theta^*(m)}^0)^2}{H_{\theta^*(m)}^0} \Big ]\quad \mbox{by Jensen Inequality.}
\end{eqnarray*}
Since $ x-\log(x)- 1 > 0$ for any $x>0,~x\neq 1$ and $ x-\log(x)- 1= 0$ for $x=1$, we deduce that
\begin{itemize}
	\item If $f_{\theta^*}^0 \ne f_{\theta^*(m)}^0$ then $\E \Big [  \frac{ (f_{\theta^*}^0-f_{\theta^*(m)}^0)^2}{H_{\theta^*(m)}^0} \Big ]>0$ and $A(m)>0.$
	\item Otherwise, then
	\begin{equation*}
	A(m) =\E \Big [ \frac{H_{\theta^*}^0}{H_{\theta^*(m)}^0}- \log\Big (\frac{H_{\theta^*}^0}{H_{\theta^*(m)}^0}\Big )-1  \Big ],
	\end{equation*}
	and from assumption \textbf{A1}, since $\theta^*\notin \Theta(m)$ and $f_{\theta^*}^0=f_{\theta^*(m)}^0$, we necessarily have $H_{\theta^*}^0\ne H_{\theta^*(m)}^0$ so that $\frac{H_{\theta^*}^0}{H_{\theta^*(m)}^0}\ne 1.$ Then $A(m)>0$.
\end{itemize}
As a consequence,

\begin{align*}
\frac {C(m)- C(m^*)}{n}=  A(m)+\frac{2\,c\, \kappa_n}{n}(D_m-D_{m^*})+ o_{a.s}(1).
\end{align*}
That establishes \eqref{eq:ous} by virtue of \eqref{eq:pena} and as all the considered models are finite dimensional.
\qed


~\\
In the sequel, we state and prove several lemmas to which we referred to when proving above main results.

\subsubsection{Proof of Theorem \ref{theo:cons2}}

The proof follows mutatis mutandis from the Theorem \ref{theo:cons}'s proof after replacing line by line $L_n$, $\widehat{\theta}(m)$ and the criterion $C$ by their equivalent $\widehat{L}_n$, $\widetilde{\theta}(m)$ and $\widehat{C}$ respectively and applying Lemma \ref{lem:pap11} instead of Proposition \ref{dodo}.
\qed

\subsubsection{Proof of Proposition \ref{dodo}}
\begin{proof}
	
	Applying a second order Taylor expansion of $L_n$ around $\widehat{\theta}(m)$ for $n$ sufficiently large such that $\overline{\theta}(m) \in \Theta(m)$ which are between $\theta^*_m:=(\theta^*,0,\ldots,0)^\top$ and $\widehat{\theta}(m)$  yields (as $\partial L_n(\widehat{\theta}(m))=0$ since $\widehat{\theta}(m)$ is a local extremum):
	\begin{eqnarray*} \label{dtheta}
		L_n(\theta^*)-	L_n(\widehat{\theta}(m))&=&L_n(\theta^*_m)-	L_n(\widehat{\theta}(m))\\
		&	=&\frac{1}{2}\big(\widehat{\theta}(m)-\theta^*_m\big)^\top\, \frac{\partial^2 L_n(\overline{\theta}(m))}{\partial \theta^2} \, \big(\widehat{\theta}(m)-\theta^*_m\big):=I_1(m).
	\end{eqnarray*}
	\noindent From the mean value Theorem, and for large $n$, there exists $\overline \theta _{m,i}$ between $(\theta^*_m)_i$ and $(\widehat \theta(m))_i$ such that, $1\le i \le D_m$:
	\begin{equation}\label{eq:mean}
	0=\frac{\partial L_n ( \widehat \theta _m)}{\partial \theta_i}= \frac{\partial L_n (\theta^*_m)}{\partial \theta_i}+\frac{\partial^2 L_n (\overline \theta _{m,i})}{\partial \theta \partial \theta_i}(\widehat \theta _m-\theta^*_m)
	\end{equation}
	\noindent Also, using Lemma 4 of \cite{barW} and continuous mapping Theorem, we deduce that:
	\begin{equation}
	F_n:=-\Big(\frac{2}{n}\frac{\partial^2 L_n(\overline \theta _{m,i})}{\partial \theta \partial \theta_i}\Big)_{1\le i\le D_m} \limiteasn F(\theta^*_m)=\E\Big[\frac{\partial^2 q_0(\theta_m^*)}{\partial \theta^2} \Big].
	\label{eq:oli}
	\end{equation}
	On the other hand, under \textbf{A2} condition, $F(\theta^*_m)$ is an invertible matrix and there exists $n$ sufficiently large such that $F_n$ is invertible. Therefore, from \eqref{eq:mean}, it follows

	\begin{eqnarray*}
		\frac{I_1(m)}{2\,\log \log n}&=&\frac{1}{4\, \log \log n} \big(\widehat{\theta}(m)-\theta^*_m\big)^\top\, \frac{\partial^2 L_n(\overline{\theta}(m))}{\partial \theta^2} \, \big(\widehat{\theta}(m)-\theta^*_m\big)\\
		&=& \frac{1}{4\,\log \log n} \, \Big( \frac{\partial L_n (\theta^*_m )}{\partial \theta}\Big)^\top \,\Big(-\frac{2}{n}\, F_n^{-1}\Big) \, \frac{\partial^2 L_n(\overline{\theta}(m))}{\partial \theta^2} \,\Big(-\frac{2}{n} \, F_n^{-1}\Big)\,  \,\frac{\partial L_n (\theta^*_m)}{\partial \theta} \,    \\
		&=&- \frac{1}{2n \,\log \log n}\Big( \frac{\partial L_n (\theta^*_m )}{\partial \theta}\Big)^\top\, \times \, F(\theta^*_m)^{-1}\, \times  \frac{\partial L_n (\theta^*_m )}{\partial \theta}\,\big(1+o(1)\big) \quad \mbox{a.s.}\\
	\end{eqnarray*}
	The next step of the proof consists in handling the quadratic form $\Big( \frac{\partial L_n (\theta^*_m )}{\partial \theta}\Big)^\top\, \times \, F(\theta^*_m)^{-1}\, \times \Big( \frac{\partial L_n (\theta^*_m )}{\partial \theta}\Big)$ by applying the law of iterated logarithm (LIL).
	~\\We claim that
	\begin{equation}\label{sto3}
	\limsup\limits_{n \rightarrow \infty}  \frac{1}{\sqrt{2n \,\log \log n}}\,2\,G(\theta^*_m)^{-1/2}  \frac{\partial L_n (\theta^*_m )}{\partial \theta}=\big(1, \ldots, 1\big)^\top.
	\end{equation}
	\textbf{Proof of the claim:} First, since the covariance matrix of $2\, \frac{\partial L_n (\theta^*_m )}{\partial \theta}$ is  $G(\theta^*_m)$, it follows that the covariance matrix of the vector $Z_m:=2\,G(\theta^*_m)^{-1/2}\,\frac{\partial L_n (\theta^*_m )}{\partial \theta}$ is the $D_m\times D_m$ identity matrix. Moreover, as 
	$$\frac{\partial L_n (\theta^*_m )}{\partial \theta}= -\frac{1}{2}\, \sum_{t=1}^n  \frac{\partial q_t (\theta^*_m )}{\partial \theta}=-\frac{1}{2}\, \sum_{t=1}^n  \frac{\partial q_t (\theta^*)}{\partial \theta},$$
	where
	\begin{equation}\label{sto1}
	\E\Big[\frac{\partial q_t (\theta^*)}{\partial \theta} \Big| \sigma(X_{t-1},X_{t-2},\ldots,X_{1})  \Big]=0   
	\end{equation}
	and
	\begin{equation}\label{sto2}
	\E\Big[\Big(\frac{\partial q_1 (\theta^*)}{\partial \theta_i}\Big)^2 \Big] < \infty
	\end{equation}
	hold from \cite{barW} under \textbf{A4}. Finally, one can see that the $i^{th}$ element of $Z_m$ can be rewritten as
	$$\sum_{j=1}^{D_m}  \big(2\,G(\theta^*_m)^{-1/2} \big)_{i,j}\,\frac{\partial L_n (\theta^*_m )}{\partial \theta_j}=\sum_{t=1}^n \zeta_t^i $$
	where $\zeta_t^i=-\sum_{j=1}^{D_m}  \big( G(\theta^*_m)^{-1/2} \big)_{i,j}\,\frac{\partial q_t (\theta^*)}{\partial \theta_j}$.
	By virtue of \eqref{sto1}, we have 
	$$ \E\Big[\zeta_t^i \Big| \sigma(X_{t-1},X_{t-2},\ldots,X_{1})  \Big]=0.$$
	Hence, any component of $Z_m$ verifies the LIL. That is for any $i=1,\ldots,D_m$,
	$$\limsup\limits_{n \rightarrow \infty}  \frac{1}{\sqrt{2n \,\log \log n}}\,\Bigg(2\,G(\theta^*_m)^{-1/2}  \frac{\partial L_n (\theta^*_m )}{\partial \theta}\Bigg)_i=1.$$
	This fact  concludes the proof of the claim \eqref{sto3}.
	
	Hence writing 
	\begin{multline*}
	\frac{1}{2n \,\log \log n}\frac{\partial L_n (\theta^*_m )}{\partial \theta}^\top  F(\theta^*_m)^{-1} \frac{\partial L_n (\theta^*_m )}{\partial \theta}= \\\Bigg(\frac{1}{\sqrt{2n \,\log \log n}}\,2G(\theta^*_m)^{-1/2}  \frac{\partial L_n (\theta^*_m )}{\partial \theta}\Bigg)^\top \, \frac{G(\theta^*_m)^{1/2} F(\theta^*_m)^{-1} G(\theta^*_m)^{1/2}}{4}  \Bigg( \frac{1}{\sqrt{2n \,\log \log n}}\,2G(\theta^*_m)^{-1/2} \frac{\partial L_n (\theta^*_m )}{\partial \theta}\Bigg)
	\end{multline*}
	it follows
	$$ \limsup\limits_{n \rightarrow \infty} \, \frac{L_n\big(\widehat{\theta}(m)\big)-L_n(\theta^*_m)}{2\,\log \log n} =\textnormal{\textbf{1}}_m^\top\,\Sigma_{\theta^*_m} \,\textnormal{\textbf{1}}_m.
	$$
\end{proof}

\subsubsection{Proof of Proposition \ref{dodo2}}
It is sufficient to show that
\begin{equation} \label{mouss}
\Sigma_{\theta^*_m}:=G(\theta^*_m)^{1/2} F(\theta^*_m)^{-1} G(\theta^*_m)^{1/2}= \left (
\begin{array}{cc}
2\, \ii_{D_m^1,D_m^1}  &  \oo_{D_m^1,D_m^2} \\
\oo_{D_m^2,D_m^1}  &  (\mu_4-1)\,\ii_{D_m^2,D_m^2}
\end{array}
\right ),
\end{equation}
where $\ii$ is the identity matrix and $\oo$ the null matrix.
From \cite{barW}, we have for $m^*\subset m$ and $i,j \in m$:
\begin{eqnarray}
\nonumber  \big (G(\theta^*_m)\big )_{i,j}&\hspace{-3mm}=&\hspace{-3mm}\E\Big[\frac{\partial q_0(\theta_m^*)}{\partial \theta_i} \frac{\partial q_0(\theta^*_m)}{\partial \theta_j} \Big]=\E\Big[4 (H_{\theta_m^*}^0)^{-1} \, \frac{\partial f_{\theta_m^*}^0}{\partial \theta_i}\,\frac{\partial f_{\theta_m^*}^0}{\partial \theta_j} + (\mu_4 -1)\, (H_{\theta_m^*}^0)^{-2} \, \frac{\partial H_{\theta_m^*}^0}{\partial \theta_i}\,\frac{\partial H_{\theta_m^*}^0}{\partial \theta_j}\Big] \\
\big (F(\theta^*_m)\big )_{i,j}&\hspace{-3mm}=&\hspace{-3mm}\E\Big[\frac{\partial^2 q_0(\theta_m^*)}{\partial \theta_i \partial \theta_j}\Big]=\E\Big[2 (H_{\theta_m^*}^0)^{-1} \, \frac{\partial f_{\theta_m^*}^0}{\partial \theta_i}\,\frac{\partial f_{\theta_m^*}^0}{\partial \theta_j} + (H_{\theta_m^*}^0)^{-2} \, \frac{\partial H_{\theta_m^*}^0}{\partial \theta_i}\,\frac{\partial H_{\theta_m^*}^0}{\partial \theta_j}\Big]
\label{eq:m1}
\end{eqnarray}

\noindent 1/ If $\mu_4=3$, then $G(\theta^*_m)=2\,F(\theta^*_m)$ and the result is straightforward.
\\
~\\
\noindent 2/  For the second configuration, from \cite{bardet2021}, we have
\begin{multline*}
G(\theta^*_m)= \left (
\begin{array}{cc}
(G(\theta^*_m)\big )_{1\le i,j \le D_m^{1}} &  \oo_{D_m^1,D_m^2} \\
\oo_{D_m^2,D_m^1}  &  (\mu_4-1)\,(G(\theta^*_m)\big )_{1\le i,j \le D_m^{2}}
\end{array}
\right )  ~
\\\mbox{and} \quad
G(\theta^*_m)F(\theta^*_m)^{-1}= \left (
\begin{array}{cc}
2\, \ii_{D_m^1, D_m^{1}} &  \oo_{D_m^1,D_m^2} \\
\oo_{D_m^2,D_m^1}  & (\mu_4-1)\,\ii_{D_m^2,D_m^2}
\end{array}
\right ).
\end{multline*}

As a covariance matrix,  $G(\theta^*_m)$ is positive definite. 
Therefore  the square root $G(\theta^*_m)^{1/2}$ of $G(\theta^*_m)$ is unique and  blocks diagonal.
Thus,
\begin{eqnarray*}
	\Sigma_{\theta^*_m}&=&G(\theta^*_m)^{-1/2}\,\Big( \,G(\theta^*_m) F(\theta^*_m)^{-1}\Big) G(\theta^*_m)^{1/2} \\
	&=&G(\theta^*_m) F(\theta^*_m)^{-1},
\end{eqnarray*}
which gives \eqref{mouss}.
\qed

\subsection{Technical Lemmas}

\begin{lemma} \label{lem11}
	Let $X \in \mathcal{H}(M_{\theta^*},f_{\theta^*})$ (or $\widetilde{\mathcal{H}}(H_{\theta^*})$) and $\Theta \subseteq \Theta(r)$  (or $\Theta \subseteq \widetilde \Theta(r)$) with {$r\ge 4$}. Assume that assumption \textbf{A3} holds. Then for $i=0,1,2$, it holds
	\begin{equation}
	\frac{1}{\log \log n} \,\Big \|\dfrac{\partial^{(i)}\widehat{L}_n(\theta)}{\partial \theta^i}-\frac{\partial^{(i)} L_n(\theta)}{\partial \theta^i}\Big \|_{\Theta}\limiteasn 0.
	\label{eq:22}
	\end{equation}
\end{lemma}
\begin{proof}
	This Lemma has already been proved in \cite{bar2019} in a more general framework. Let us prove the result for $i=0$. Other cases can be deduced by using a similar reasoning. \\
	~\\
	We have for any $\theta \in \Theta$,
	$|\widehat{L}_n(\theta)-L_n(\theta)|\le \sum_{t=1}^n|\widehat{q}_t(\theta)-q_t(\theta)| $. Then,
	\[
	\frac{1}{\log \log n} \, \big \|\widehat{L}_n(\theta)-L_n(\theta) \big \|_{\Theta}\le \frac{1}{\log \log n} \, \sum_{t=1}^n\big\|\widehat{q}_t(\theta)-q_t(\theta)\big\|_{\Theta}.\]
	By Corollary 1 of \cite{kounias}, \eqref{eq:22} is established when:
	\begin{equation}
	\sum_{k\ge 1}\frac{1}{\log \log k}\, \mathbb{E}\big\|\widehat{q}_k(\theta)-q_k(\theta)\big\|_{\Theta} < \infty.
	\end{equation}
	From \cite{barW} and \cite{bar2019}, there exists a constant $C$ such that

	\noindent 1/ If $X \in {\mathcal H}(M_\theta,f_\theta)$, we deduce

	\begin{equation}
	\mathbb{E} \big [ \|(\widehat{q}_{t}(\theta)-q_{t}(\theta)\|_{\Theta} \big ] \le C \,\Big(\sum_{j \ge t} \alpha_j (f_\theta,\Theta)+\sum_{j \ge t} \alpha_j (M_\theta,\Theta)\Big).
	\label{e:fond}
	\end{equation}
	Hence,
	\[
	\sum_{k\ge 1}\frac{1}{\log \log k} \mathbb{E}\big[\|\widehat{q}_k(\theta)-q_k(\theta)\|_{\Theta} \big] \le C \, \sum_{k\ge 1}\frac{1}{\log \log k}\Big(\sum_{j \ge k} \alpha_j (f_\theta,\Theta)+\alpha_j (M_\theta,\Theta)\Big),
	\]
	which is finite  by definition of the class $\mathcal{H}(M_{\theta},f_{\theta})$, and this achieves the proof. \\
	~ \\
	\noindent 2/ If $X \in \widetilde {{\mathcal H}}(\widetilde H_\theta)$,
	%
	\begin{equation}
	\mathbb{E} \big [ \|(\widehat{q}_{t}(\theta)-q_{t}(\theta)\|_{\Theta} \big ] \le C \,\Big(\sum_{j \ge t} \alpha_j (H_\theta,\Theta)\Big).
	\label{e:fond2}
	\end{equation}
	This fact along with Corollary 1 of \cite{kounias} enable us to conclude the proof in this case.
\end{proof}

\begin{lemma} \label{lem:pap11}
	Under the assumptions of Theorem \ref{theo:cons2}, for any model $m\in\mathcal{M}$ with  $\theta^*\in \interior{\widetilde{\Theta(m)}}$, it holds 
	\begin{eqnarray}\label{eq:kami2}
	\limsup\limits_{n \rightarrow \infty} \Bigg(\frac{\widehat{L}_n\big(\widetilde{\theta}(m)\big)-\widehat{L}_n(\theta^*)}{2\,\log \log n}   \Bigg) & = & c_{min}(m) \quad a.s.
	\end{eqnarray}
\end{lemma}
\begin{proof}
	Applying a second order Taylor expansion of $\widetilde{L}_n$ around $\widetilde{\theta}(m)$ for $n$ sufficiently large such that $\overline{\theta}(m) \in \Theta(m)$ which are between $\theta^*_m:=(\theta^*,0,\ldots,0)^\top$ and $\widetilde{\theta}(m)$  yields (as $\partial \widehat{L}_n(\widetilde{\theta}(m))=0$ since $\widetilde{\theta}(m)$ is a local extremum):
	\begin{eqnarray*} 
		\widehat{L}_n(\theta^*)-\widehat{L}_n\big(\widetilde{\theta}(m)\big)&=&\widehat{L}_n(\theta^*_m)-\widehat{L}_n\big(\widetilde{\theta}(m)\big)\\
		&	=&\frac{1}{2}\big(\widetilde{\theta}(m)-\theta^*_m\big)^\top\, \frac{\partial^2 \widehat{L}_n(\overline{\theta}(m))}{\partial \theta^2} \, \big(\widetilde{\theta}(m)-\theta^*_m\big)\\
		&:=&  I_2(m).
	\end{eqnarray*}
	But $I_2(m)$ can be rewritten as
	\begin{multline*}
	\frac{I_2(m)}{2\,\log \log n}=\frac{1}{4}\,\big(\widetilde{\theta}(m)-\theta^*_m\big)^\top\,\frac{1}{\log \log n}\Bigg[ \frac{\partial^2 \widehat{L}_n(\overline{\theta}(m))}{\partial \theta^2} -\frac{\partial^2 L_n(\overline{\theta}(m))}{\partial \theta^2}\Bigg]\, \big(\widetilde{\theta}(m)-\theta^*_m\big)\\+\frac{1}{4}\,\frac{1}{\log \log n}\big(\widetilde{\theta}(m)-\theta^*_m\big)^\top\, \frac{\partial^2 L_n(\overline{\theta}(m))}{\partial \theta^2} \, \big(\widetilde{\theta}(m)-\theta^*_m\big)\\
	=:I_{21}(m)+I_{22}(m).
	\end{multline*}
	First, as $\widetilde \theta(m)\limiteasn \theta^*_m $ along side with  Lemma \ref{lem11}, it follows
	\begin{equation}
	I_{21}(m) \limiteasn 0.
	\end{equation}
	Writing a counterpart of \eqref{eq:mean} using the quasi-functions, we have
	\begin{eqnarray*}
		\nonumber \widetilde \theta(m)-\theta^*_m&=&  \Big(\partial^2 \widehat{L}_n(\overline{\theta}(m)) \Big)^{-1} \frac{\partial \widehat{L}_n (\theta^*_m)}{\partial \theta}   \\
	\end{eqnarray*}
	Hence $I_{22}(m)$ can be rewritten as
	\begin{eqnarray}
	\nonumber I_{22}(m) &=& \frac{1}{4\,\log \log n}\big(\widetilde{\theta}(m)-\theta^*_m\big)^\top\, \frac{\partial^2 L_n(\overline{\theta}(m))}{\partial \theta^2} \, \big(\widetilde{\theta}(m)-\theta^*_m\big) \\
	\nonumber &=&\frac{1}{4\,\log \log n} \,\Big( \frac{\partial \widehat{L}_n (\theta^*_m )}{\partial \theta}\Big)^\top  \Big(\partial^2 \widehat{L}_n(\overline{\theta}(m)) \Big)^{-1} \, \Big(\partial^2 L_n(\overline{\theta}(m))\Big) \, \Big(\partial^2 \widehat{L}_n(\overline{\theta}(m)) \Big)^{-1} \,  \,\frac{\partial \widehat{L}_n (\theta^*_m)}{\partial \theta}\\
	&=& - \frac{1}{2\,n \,\log \log n}\Big( \frac{\partial L_n (\theta^*_m )}{\partial \theta}\Big)^\top\, \times \, F(\theta^*_m)^{-1}\, \times  \frac{\partial L_n (\theta^*_m )}{\partial \theta}\big(1+o(1)\big) \quad \mbox{a.s.}
	\end{eqnarray}
	since from \eqref{eq:oli},  it holds
	$$ -\frac{n}{2}\, \Big(\partial^2 L_n(\overline{\theta}(m)) \Big)^{-1} \limiteasn F(\theta^*_m)^{-1}$$
	and along with Lemma \ref{lem11}, it also holds
	$$-\frac{n}{2} \log \log  n \, \Big(\partial^2 \widehat{L}_n(\overline{\theta}(m)) \Big)^{-1} \limiteasn F(\theta^*_m)^{-1}.$$
	
	As a consequence, the chain of following equalities holds a.s.
	\begin{eqnarray*}
		\limsup\limits_{n \rightarrow \infty} \Bigg(\frac{\widehat{L}_n\big(\widetilde{\theta}(m)\big)-\widehat{L}_n(\theta^*)}{2\,\log \log n}   \Bigg)&=&\limsup\limits_{n \rightarrow \infty} \Bigg( - \frac{1}{2\,n \,\log \log n}\Big( \frac{\partial L_n (\theta^*_m )}{\partial \theta}\Big)^\top\, \times \, F(\theta^*_m)^{-1}\, \times  \frac{\partial L_n (\theta^*_m )}{\partial \theta}\Bigg)\\
		&=&c_{min}(m)
	\end{eqnarray*}
	That ends the proof of \eqref{eq:kami2}.
\end{proof}


\section{Acknowledgements}
The author thanks Jean-Marc BARDET for proofreads and helpful discussions.

\newpage
\small{
\bibliographystyle{abbrv}
\bibliography{biblio7}}
\end{document}